\documentclass[12pt,a4paper,reqno]{amsart}

\usepackage[margin = 1in]{geometry}
\usepackage{a4wide}
\usepackage{amssymb,amsmath,amsthm,mathrsfs}
\usepackage{graphicx}
\usepackage{verbatim}
\usepackage{caption}
\usepackage{float}
\usepackage{wrapfig}
\usepackage{colortbl}
\usepackage{cite}
\usepackage{enumerate}
\usepackage{upref}
\usepackage{tikz}
\usepackage{pgfplots}
\pgfplotsset{compat=1.11}
\usepgfplotslibrary{fillbetween}
\usetikzlibrary{intersections}
\usetikzlibrary{arrows,backgrounds,patterns}
\pgfdeclarelayer{bkgr}
\pgfdeclarelayer{bkgrd}
\pgfsetlayers{bkgrd,bkgr,main}
\usetikzlibrary{arrows}
\usetikzlibrary{decorations.pathmorphing}
\usetikzlibrary{decorations.pathreplacing}
\usepackage[T2A,T1]{fontenc}
\DeclareSymbolFont{cyrillic}{T2A}{cmr}{m}{n}
\DeclareMathSymbol{\D}{\mathalpha}{cyrillic}{196}

\usepackage{bbm}
\usepackage[bookmarks=false]{hyperref}
\theoremstyle{plain}
\newtheorem{theorem}{Theorem}[section]
\newtheorem{lemma}[theorem]{Lemma}

\theoremstyle{definition}
\newtheorem{definition}{Definition}

\theoremstyle{remark}
\newtheorem{remark}[theorem]{Remark}

\usepackage{enumitem}
\makeatletter
\def\namedlabel#1#2{\begingroup
   #2%
 \def\@currentlabel{#2}%
   \phantomsection\label{#1}\endgroup
}

\usepackage[T2A,T1]{fontenc}
\makeatletter
\newcommand*{\math@version@bold}{bold}
\DeclareMathOperator\DD{
  \textrm{%
    \usefont{T2A}{cmr}{\ifx\math@version\math@version@bold bx\else m\fi}{n}%
    \CYRD
  }%
}
\makeatother

\numberwithin{equation}{section}

\begin{document}

\title[\hfill\protect\parbox{0.975\linewidth}{Multiday Extreme Precipitaton Across Eastern Australia. }]{Modeling Multiday Extreme Precipitation Across Eastern Australia: A Dynamical Perspective.}

\author[R. Kardkasem]{Ruethaichanok Kardkasem}
\address{Ruethaichanok Kardkasem\\ School of Mathematics and Physics\\
University of Queensland\\
St Lucia\\
QLD 4067\\
USA} \email{r.kardkasem@uq.edu.au}
\urladdr{}

\author[M. Carney]{Meagan Carney$^{\ast}$}
\address{Meagan Carney\\ School of Mathematics and Physics\\
University of Queensland\\
St Lucia\\
QLD 4067\\
USA} \email{m.carney@uq.edu.au}
\urladdr{https://smp.uq.edu.au/profile/11136/meagan-carney}

\thanks{$^{\ast}$ Corresponding author.}

\begin{abstract}
The purpose of this paper is to illustrate new techniques for computing multiday extreme precipitation taken from recent theoretical advancements in extreme value theory in the framework of dynamical systems, using historical precipitation data along the eastern coast of Australia as a case study. We explore the numerical pitfalls of applying standard extreme value techniques to model multiday extremes. Then, we illustrate that our data conforms to the appropriate setting for the application of recently derived extreme value distributions for runs of extremes in the dynamical framework and adapt these to the non-stationary setting. Finally, we use these distributions to make more informed predictions on the return times and magnitudes of consecutive daily extreme precipitation and find changes in the dependence of increasing consecutive daily rainfall extremes on the Southern Oscillation Index. Although our case study is focused on extreme precipitation across eastern Australia, we emphasize that these techniques can be used to model expected returns and magnitudes of consecutive extreme precipitation events across many locations.
\end{abstract}

\maketitle

\section{Introduction and Background}
\subsection*{Summary of extremes in rainfall across Australia.}\label{intro1}

In recent decades, eastern Australia has experienced flooding attributed to multiday extreme rainfall, for example, the historic 1974 floods in Brisbane, the Queensland and Victoria floods occurring through 2010 and 2011, and the 2022 flooding across New South Wales and Queensland. These catastrophic flooding events are considered some of the most substantial economic risks to the country, with the Insurance Council of Australia reporting that the 2022 floods were the costliest in Australia's history, resulting in an estimated \$3.35 billion in insured losses \cite{ICRR}. Accurately estimating returns and magnitudes of multiday extreme rainfall can inform infrastructure, city planning, and agriculture to help combat expected future loss.

Precipitation extremes are considered primary triggers for flooding alongside other factors, such as climate, topography, and soil type of different catchment patterns. Precipitation extremes along eastern Australia are driven by effects of El Nino-Southern Oscillation (ENSO) \cite{csiro2015climate, king2013limited, van20112011}. La Niña events are often blamed for higher multiday rainfall, which tends to be attributed to devastating floods; however the effects of ENSO vary significantly over geographic location \cite{bureau2022special, national2011frequent, van20112011}. Several studies have analyzed extreme precipitation in Australia. Sun et al. (2021) \cite{sun2021global} conducted a study that found a decrease in one-day and five-day precipitation extremes in Australia, in contrast to the global trend of increasing precipitation extremes. However, Alexander and Arblaster (2017) \cite{alexander2017historical} discovered that global climate models show increased precipitation extremes across Australia. In another study, Hajani et al. (2017) \cite{hajani2017trends} examined the annual maximum of sub-daily to daily precipitation in New South Wales, Australia, and found that the number of stations showing significant trends decreased when the impact of climate variability modes was considered. Jakob and Walland (2016) \cite{jakob2016variability} revealed that there has been no significant change in daily precipitation extremes in Australia over the past century, but the ENSO significantly impacts the occurrence of these extremes. Laz et al. (2014) \cite{laz2014trends} analyzed extreme rainfall patterns from 1950-2010 in southeast Australia on both sub-daily and daily scales and discovered that annual maximum rainfall events lasting less than 12 hours increased, while rainfall events lasting 12-72 hours decreased. Westra and Sisson (2011) \cite{westra2011detection} reported a similar result. In particular, they analyzed the annual maxima precipitation in eastern Australia from 1965 to 2005 and reported a statistically significant increase in short-duration rainfall over this period but no trend for daily and a decreasing trend in multiday precipitation. Many of these studies have explored trends in extreme precipitation, demonstrating the nonstationarity in precipitation time-series in relation to ENSO.

A common approach for modeling extremes in a sequence is to use results from the field of Extreme Value Theory which provide distributional assumptions for the extremes of time-series under certain weak dependence conditions, provided the sequence is stationary or does not change in time or space. Consequently, the spatial and temporal nonstationarity observed in weather data introduces complexities in modeling using standard extreme value statistics; however, many studies have successfully addressed this issue by including nonstationary dependencies directly into the model to create a unified spatiotemporal model of rainfall extremes. Recent studies developed a spatiotemporal nonstationary model for extreme precipitation across regions in Australia, including large-scale modes of climate variability and topology variables \cite{yilmaz2017investigation, saunders2017spatial, westra2011detection,jayaweera2023non}. However, there are limited studies concerning consecutive days of extreme precipitation in the model. Despite \cite{jayaweera2023non, westra2011detection} applying a nonstationary generalized extreme value (GEV) distribution to model extreme precipitation across several durations from 6 minutes through to 72 hours annual maxima, the studies did not provide the goodness of fit of the model. In addition, it is not always clear how concurrent daily extreme precipitation is influenced by the commonly used indices for monitoring extreme precipitation, such as annual maximum consecutive 5-day precipitation (Rx5day), annual total precipitation from days greater than the 99th percentile (R99p), or maximum number of consecutive days with daily precipitation greater than 1 mm (CWD) (see e.g., \cite{alexander2017historical, hajani2017trends, laz2014trends, sun2021global, westra2011detection, yilmaz2017investigation}). For example, it is known that increasing the number of consecutive days one considers for the model will reduce the data available for fitting the model by the same factor. Consequently, the reliability of the statistical model decreases significantly as the investigated number of consecutive days of extreme rainfall increases.

\subsection*{Dynamical systems and recent advancements in modeling extreme weather.}
Recent advancements in the field of dynamical systems have provided insights into key factors for more accurately modeling the consecutive nature of extreme weather events. To explain the impact of these results, we begin by introducing some background on the relationship between extreme value theory for dynamical systems and their role in deriving extreme value statistics for climate applications. 

Beginning in the late 1970s with the work of Leadbetter and Chernick, \cite{leadbetter1983extremes, chernick1991calculating} an extreme value law guaranteeing the existence of an extreme value distribution was proven in the limit for weakly dependent sequences satisfying a mixing and recurrence condition called $D(u_n)$ and $D^k(u_n)$, respectively \cite{chernick1991calculating}. Moreover, Leadbetter \cite{leadbetter1983extremes} showed that under these conditions the resulting extreme value law was similar to the one derived for independent and identically distributed random variables apart from the introduction of a new parameter, $\theta\le 1$ called the extremal index, which characterizes the short-term dependence observed in the weakly-dependent sequence. The authors in \cite{hsing1993extremal} showed under general conditions that the returns of a the sequence to an extreme or rare event can be understood as some (compound) Poisson process and $1/\theta$ represents the expected cluster size of extremes. The formal definition of the Extreme Value Law is provided below.
\begin{definition}
Extreme Value Law (EVL) Let $(u_n)$ be a sequence of constants defined by the requirement that $\lim_{n\rightarrow\infty}n P(X_1>u_n) = \tau$, and suppose $(X_n)$ is stationary and satisfies $D(u_n)$ and $D^k(u_n)$, then, 
\[
\lim_{n\rightarrow\infty} P(M_n\le u_n )=e^{-\theta\tau}
\]
where the sequence $M_n=\max\{X_1,…,X_n\}$ and parameter $\theta\le 1$ is the extremal index.
\end{definition}

For numerical applications of data modeling, the sequence $(u_n)$ can be understood as some sequence of increasing thresholds under which a frequency can be calculated, $u_n= ya_n+b_n$ where $(a_n)$ and $(b_n)$ are normalizing sequences with $(a_n )>0~\forall n$. For some fixed $n$, the sequence $(a_n)$ and $(b_n)$ take fixed values in $\mathbb{R}$ which can be understood as some scale $\sigma$ and location $\mu$ parameter. Under this framework, we can rewrite the EVL in the following way,  
$$P(a_n (M_n-b_n )\le y)\rightarrow G(y)$$ 
where $G(y)$ is the limiting extremal distribution of the suitably normalized maxima and takes the form of one of three types of distributions. 

A Type 1 (Gumbel) distribution has an upper exponentially decaying tail; Type 2 (Fr\'{e}chet) corresponds to an upper polynomial decaying tail; and Type 3 (Weibull) corresponds to a bounded tail. 

Results for i.i.d. sequences \cite{fisher1928limiting}, later extended for stationary sequences with weak dependence \cite{leadbetter1983extremes}, guarantee that all three extreme value distributions may be combined into the following three parameter distribution function,
\begin{definition} 
Generalized Extreme Value Distribution (GEV)
$$P(M_n\le z) \rightarrow G(z) = \exp\bigg(-\bigg[1+\xi \bigg(\frac{z-\mu}{\sigma}\bigg)\bigg]^{1/\xi}\bigg)$$
where $\mu$ is the location parameter, $\sigma$ is the scale parameter, and $\xi$ is the shape parameter which determines the tail-behavior or type of distribution \cite{Coles2001}. 
\end{definition}
One key point that should be emphasized here is the relationship between the GEV for i.i.d. sequences and the GEV for weakly dependent, stationary sequences with extremal index $\theta<1$.

\begin{theorem}[Coles \cite{Coles2001}] Given any weakly dependent sequence with extremal index $\theta$ and GEV, $G_2(\mu^*,\sigma^*,\xi)$, there exists a GEV, $G(\mu,\sigma,\xi)$ for a corresponding independent sequence such that,
\[
G_2(\mu^*,\sigma^*,\xi) = G^{\theta}(\mu,\sigma,\xi)
\]
where the extremal parameters of $G_2$ and $G$ are related by,
\[
\xi = \xi,\quad \mu^* = \mu-\frac{\sigma}{\xi}(1-\theta^{\xi}), \quad \sigma^* = \sigma\theta^{\xi}.
\]

\end{theorem}

Provided the data meet the assumptions of stationarity and weak dependence, the limiting distribution of $(M_n)$ is known, and creating a numerical model for $(M_n)$ is then reduced to applying an appropriate statistical fitting method that results in the most accurate numerical approximation of the parameters $\xi$, $\mu$, and $\sigma$. Since traditional statistical fitting methods require a frequency plot of the data to estimate parameters, the block maxima method is typically employed to generate such a frequency plot in which the data is blocked into pieces of length $m$ (long enough to guarantee convergence to the GEV), the maxima are calculated over each block, and this sequence of maxima is used to approximate the parameters of the corresponding GEV. 


Recent advancements in understanding and modeling consecutive climate extremes have centered around redefining the sequence $(X_n)$ as an average over some window and interpreting returns of extreme highs or lows of this average as indicative of several high or low values, respectively, occurring within the window. However, it was only recently that the extreme value theory and dynamics communities began investigating the numerical accuracy of these models \cite{ragone} with the first expectations on the limiting distribution type established this year \cite{CHNT}. For ease of interpretation, we focus on results for extremes of the minimum over a window; however, we recognize the importance of the windowed average in climate applications and remark that expectations on the limiting distribution for certain windowed averages have also been established \cite{CHNT}. We now restate for our context the dynamical lemma which will be central to our investigation on multiday precipitation extremes across eastern Australia.
\begin{lemma}[Lemma 4.1 \cite{CHNT}]
Suppose $M_n = \max\{X_1,\dots,X_n\}$ follow a GEV, $G(\xi_1,\mu_1,\sigma_1)$ of Fr\'{e}chet type. Define $Y_j = \min\{X_j,\dots,X_{j+k-1}\}$ for $j = 1,\dots n$ to be the minimum value over a time window of length $k$ and $B_n = \max\{Y_1,\dots Y_n\}$, so that if $B_i\ge z$ for $z\in\mathbb{R}$ then every value over the window of length $k$, ensures $X_i\ge z,\dots,X_{i+k-1}\ge z$. Suppose $(X_n)$ and $(Y_n)$ are stationary, weakly-dependent sequences, satisfying $D(u_n)$ and $D^k(u_n)$, with extremal index $\theta_1$ and $\theta_2$, respectively. Then the maximal process $(B_n)$ follows a GEV, $G(\xi_2, \mu_2,\sigma_2)$ of Fr\'{e}chet type with,
\[
\xi_2 = \xi_1, \quad \mu_2 = g(k,\chi)\mu_1\bigg(\frac{\theta_2}{\theta_1}\bigg)^{\xi_1}, \quad \sigma_2 = g(k,\chi)\sigma_1\bigg(\frac{\theta_2}{\theta_1}\bigg)^{\xi_1}
\]
for some unknown function $g(k,\chi)$ depending on the window size $k$ and the system $\chi$. 
\end{lemma}

Introducing functionals, such as the windowed minimum or average, reduces the amount of data available for fitting by a factor of the window size. Hence, for increasing window sizes, even for window sizes as low as $k = 2$ depending on the available pool of data, the numerical error of the extremal parameter estimates significantly increases. To our knowledge, this manuscript is the first to use techniques that reduce numerical error when using functionals to estimate consecutive extreme weather events. The fundamental importance of the previous lemma is that it proves max-stability (e.g. stability of the shape parameter) for the maximal process $(B_n)$. As a result, this information provides us with greater control over the numerical error of our model by forcing us to only choose models of $(B_n)$ with shape parameters \textit{close} to the true shape estimate from $(M_n)$. In the sections to follow, we illustrate where this information can be used in practice to obtain more robust modeling results on the basis of model stability. 

\section{Data and Methodology}
\subsection*{Description of the data.}
 Daily precipitation data for eastern Australia from 1984 to 2021 was obtained from the Australian Bureau of Meteorology. The study area includes New South Wales (NSW), Victoria (VIC), Queensland (QLD), and Tasmania (TAS).  The dataset consists of 454 stations that have long-term records of precipitation observations with less than 1\% missing values each year. The full list of stations is available in the Supplementary Material.


To assess the impact of ENSO on the distribution of extreme events, we used the monthly Southern Oscillation Index (SOI) to provide information about the development and intensity of El-Niño or La-Niña events in the Pacific Ocean obtained from the NOAA Climate Prediction Center.
\subsection*{Investigating extreme precipitation clustering via the extremal index.}
Formally, the extremal index, $\theta$ in the definition of an EVL, of a weakly-dependent sequence characterizes the short-term dependence observed in the recurrence condition $D^k(u_n)$. Under the assumption of a compound-Poisson limiting distribution for returns, the authors in \cite{hsing1993extremal} have shown that $1/\theta$ can be interpreted as the expected cluster size of exceedances. That is, when one exceedance occurs, one would expect that $1/\theta$ exceedances on average occur in a short window. There is extensive literature on the interpretation and approximation of the extremal index. We refer the reader to \cite{ExtremesBook} for a nice summary of the current methods used for numerical approximation of the extremal index. 

Given its proven numerical robustness, we have chosen to approximate the extremal index via the Ferro-Segers estimate \cite{ferro2003inference}. The estimate is well-accepted in the extremal community and, for replicability of results, is performed through a function that is built into many extremal platforms (e.g. MATLAB, R, Python). The Ferro-Segers estimate assumes a limiting compound-Poisson distribution for returns in our daily precipitation time-series and approximates the extremal index through a maximum likelihood calculation of the mean intensity of the compound-Poisson distribution defined by,
\[
\theta = 2\frac{\sum_{j = 1}^{N}(T_j)^2}{N\sum_{j=1}^N(T_j)(T_j-1)}
\]
where $T_j$ is the wait time associated to the exceedance $Z_j>u$ for $j = 1,\dots N$ over some high (fixed) threshold, often taken as the 95\% or 99\% quantile of the data $(X_i)$ for $i = 1,\dots n$. 

       \begin{wrapfigure}{l}{0.35\textwidth}
	 \begin{center}
	   \vspace{-25pt}
            \includegraphics[width=0.38\textwidth]{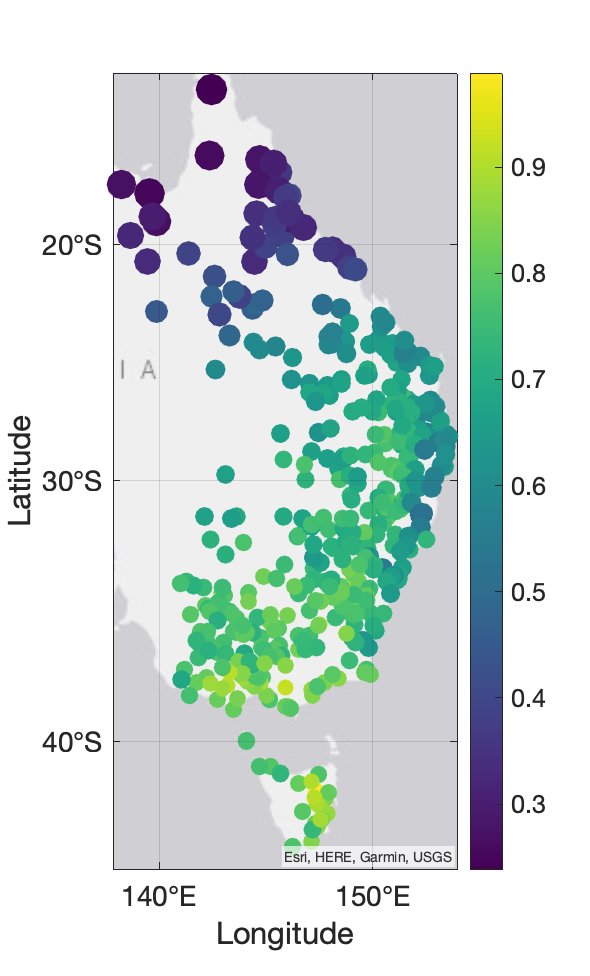}
            \vspace{-30pt}
            \end{center}
            \caption{Ferro-Segers extremal index estimate for daily rainfall from stations across eastern Australia.}
            \label{fig:ei}
             \vspace{-25pt}
        \end{wrapfigure}

Immediately obvious is the spatial structure of the extremal index from stations across eastern Australia (Figure \ref{fig:ei}). We note that the extremal index is smaller (e.g. cluster-size is larger) for stations in northern Queensland, and along the coast. This structure alone tells the story of large storms with multi-day rainfall moving in from northern Queensland and the coast, breaking up as they traverse south and inland. In 2022, the flooding across southern Queensland and northern New South Wales was attributed to the "2-3 consecutive days of heavy rainfall coming from a large, slow-moving storm system" \cite{qldflood2022}. To illustrate the importance of this index, note that the extremal index for Brisbane is estimated at $\theta\approx 0.45$ so that $1/\theta\approx 2.25$ giving an expectation of about 2.25 days of extreme rainfall which corresponds well to what was observed during the 2022 flooding event \cite{qldflood2022}. 

When numerically fitting the GEV using the block maxima method described in the previous section, one can only recover $G_2$ from Theorem 2, so the effects of the extremal index on the return probabilities are buried in the location and scale parameters of the distribution. Consequently, the physical interpretation of a return level for daily rainfall with block length of one-year is \textit{the amount of rainfall we expect a \textbf{single} yearly maximum to exceed in a given number of years}. This has clear consequences on the prediction of consecutive occurrences of extremes. One way of addressing this issue is to introduce functionals that preserve the consecutive property of extremes.

We focus our study on applications of the Generalized Extreme Value distribution; however, there is one other common approach for modeling extremes of a sequence. Before introducing our functional, we briefly remark on the appropriateness of this method and where interpretations for consecutive extremes fail.

\subsubsection*{A note on the appropriateness of the Generalized Pareto.}
Another common approach in extremal literature is to apply a peaks over threshold method which involves setting a \textit{high} threshold, creating a frequency distribution from over-threshold values in sequence, and fitting the Generalize Pareto distribution (GPD). However, the extremal index must be equal to 1 (e.g. $\theta = 1$) to guarantee a GPD models these over-threshold values. A common solution is to remove all clustered over-threshold values in the sequence, keeping only a single value to model the GPD, then factoring the extremal index back into the resulting return-level estimate by multiplying the corresponding return-time by $1/\theta$. The problem with this approach is that it removes all interpretation of consecutive extreme events, in favor of some average return-time.

\subsection*{A minimum functional to preserve the consecutive property.}

Following the theoretical work of \cite{CHNT}, we introduce the minimum functional (a moving minimum over window size $k$),
$$Y_j = \min\{X_j,\dots,X_{j+k-1}\}$$
on our daily time-series of precipitation $(X_n)$ for $j = 1,\dots n$. Provided our sequence $(X_n)$ is stationary, and blocks are taken \textit{long enough} so that observations are roughly independent and converge to the GEV, we may apply the block maxima method to estimate the maximal process $B_{i,m}$ from blocks of length $m$ given by, $B_{i,m} = \max\{Y_{im},\dots Y_{i(m+1)}\}$. For this investigation, we find block lengths $m$ of one-year long enough to guarantee convergence to the GEV.

\subsection*{Addressing non-stationarity in the data.} It is well-known that daily precipitation amounts are dependent on factors such as seasonal cycles, the strength of ENSO, and geospatial information. Taking blocks of length one-year removes seasonal and daily cycles; however, other factors such as ENSO and geographic location must be investigated and addressed. Many investigations have successfully addressed the issue of non-stationarity in climate data by including spatial and temporal variables into parameters (location, scale, and shape) of the GEV, and performing maximum likelihood estimation on their coefficients \cite{Coles2001,CKN, CK, CNA, saunders2017spatial}. In agreement with previous literature we find our GEV parameters are significantly dependent on SOI, coastal distance, latitude, and longitude in the following way,
\begin{align*}
\mu &= \mu_0+\mu_1(\text{SOI})+\mu_2\log(\text{cdist})+\mu_3(\text{lat})+\mu_4(\text{lon})\\
\sigma &= \sigma_0+\sigma_1(\text{SOI})+\sigma_2\log(\text{cdist})+\sigma_3(\text{lat})+\sigma_4(\text{lon})\\
\xi &= \xi
\end{align*}
We refer the reader to a nice, detailed investigation on modeling precipitation extremes across eastern Australian in the single-daily case, using the Generalized Extreme Value distribution \cite{saunders2017spatial}. Scatter plot for covariates against block maxima used to analyze these relationships can be found in Figure \ref{fig:scatter}.

\begin{remark}
Coefficients for the location and scale estimate decrease for increasing window size $k$. This is in line with the more robust theoretical result from \cite[Theorem 5.1]{CHNT} by noting that, rainfall is maximized at the center of a storm system (invariant set), hence $g(k,\chi)\theta_2/\theta_1$, decreases for increasing $k$. Estimating this relationship is beyond the scope of this investigation which is centered around applications of the proven max-stability property. However, a preliminary investigation using rainfall from Germany (in the approximately stationary setting) has been performed with encouraging results that support the idea of future numerical investigations. 
\end{remark}

Yearly non-stationarity was also examined as a potential covariate; however, we find that the inclusion of this effect does not result in a statistically significantly improved model over ENSO alone. We refer to Table \ref{tab:model} for model estimates, noting that we observe a Fr\'{e}chet distribution ($\xi>0$) for precipitation extremes across eastern Australia.

\begin{table}[h!]
\begin{center}
\begin{tabular}{ c | c c c c c c c c c c c}
window size & $\mu_0$ & $\mu_1$& $\mu_2$ & $\mu_3$ & $\mu_4$ & $\sigma_0$ & $\sigma_1$ & $\sigma_2$ & $\sigma_3$ & $\sigma_4$ & $\xi$\\
\hline \hline
$k$ & int & SOI & $\log(\text{cdist})$ & lat & lon & int & SOI & $\log(\text{cdist})$ & lat & lon & int\\
\hline
1 &	-71.52 &	1.28 &	-7.10	 &	2.00	&	1.48	&	58.01 &	0.09	&	-4.24	 &	1.07	&	0.11	&	0.16\\
2 &	-55.85 &	0.68 &	-3.88 &	0.77	&	0.80	&	2.96	&	0.09	&	-2.13	 &	0.55	&	0.23	&	0.18\\
3 &	-29.17 &	0.33 &	-2.29 &	0.28	&	0.39	&	-8.68	 &	0.15	&	-1.35	 &	0.29	&	0.20	&	0.24\\
\end{tabular}
\end{center}
\caption{Chosen nonstationary model coefficient estimates. There are three models, each corresponding to consecutive $k$ daily extreme rainfall. Standard errors can be found in Table \ref{tab:se}.}\label{tab:model} 
\end{table}

\subsection*{Controlling hidden numerical error of the fitted distribution of multiday extremes.}
\begin{wrapfigure}{r}{0.49\textwidth}
	 \begin{center}
            \includegraphics[width=0.48\textwidth]{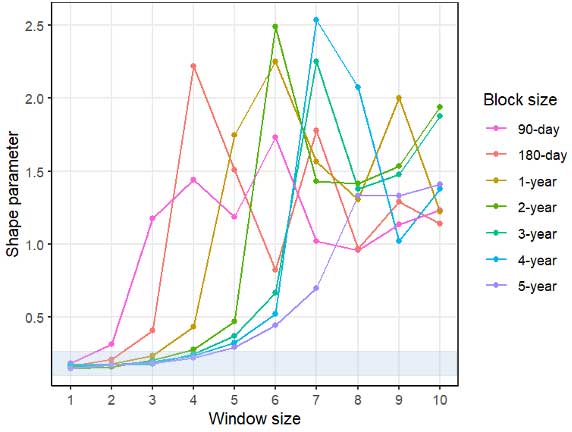}
             \vspace{-10pt}
            \end{center}
            \caption{Shape parameter estimate of the GEV for $B_{i,m}$ for increasing window size. Longer block lengths are able to capture the true shape parameter for greater window sizes.}
            \label{fig:blck}
             \vspace{-28pt}
        \end{wrapfigure}

Arguably, the most important parameter in the GEV is the shape because it provides return estimates on the tail values, or the most extreme extremes. Significant literature exists on convergence rates of the GEV, all agreeing that convergence of the shape parameter occurs at the slowest rate, forming a bottleneck of numerical accuracy. Fr\'{e}chet and Weibull distributions have significantly faster convergence $O(m^{\delta})$ $\delta>0$, over the Gumbel distribution $O(\log(m))$ \cite{smithrates}; however, introducing a windowed minimum or average functional significantly reduces the amount of tail data available to numerically fit the distribution. To our knowledge, this hidden numerical error, coming from reduced convergence rates of $(B_{i,m})$ to the GEV which is not reflected in the confidence intervals of likelihood methods, has not been considered in previous literature when forming distributions for functionals of a sequence. We illustrate this phenomenon by fitting the shape of the GEV for increasing window size over increasing block sizes in Figure \ref{fig:blck}.
        
\begin{remark}
Although we use a fixed block length for ease of interpretation and comparisons across different window sizes, a straight-forward solution to the issue of reduced convergence rates of $(B_{i,m})$ is to increase block length $m$. However, given the novelty of these results, there is currently no literature addressing optimal increases of block length.
\end{remark}
Applying the max-stability property provided by Lemma 1.2 implies that the shape parameter for the maximal process $(B_{i,m})$ should not change with increasing window size. We use this property as a benchmark for our extremal distribution fits. We find that we are able to control the shape parameter (to within 52.26\% relative error) up to window size $k = 3$, before lack of tail sampling results in a significant jump in numerical error of the shape parameter (greater than 100\% relative error). 

\subsection*{Estimating goodness-of-fit under non-stationary conditions.}
Model goodness-of-fit was examined graphically by quantile plots. They are especially effective in detecting discrepancies in the extreme upper tail of a distribution. Given an ordered sample of weakly dependent observations $x_{(1)} \leq x_{(2)} \leq ... \leq x_{(n)}$ from a population with estimated distribution function $\hat{F}$, a quantile plot consists of the points
        \begin{equation}
            \biggl\{ \hat{F}^{-1}\left( \frac{i}{n+1} \right), x_{i}; \hspace{1em} i=1,...,n \biggl\}
        \end{equation}
    If $\hat{F}$ is a reasonable estimate of $F$ or the model is a good fit, the plot should be roughly a straight line of unit slope through the origin. 
    
In the presence of non-stationary, it is necessary to use a standardized version of the data by transforming the data to a stationary random variable \cite{Coles2001}. On the basis of a non-stationary GEV distribution $Z_t \sim G(x_t; \mu_t, \sigma_t, \xi_t)$, where $t$ indicates dependence on the model covariate(s), the standardized variable $\Tilde{Z}_t$ is defined by
            \begin{equation}
                \Tilde{Z}_t = \frac{1}{\xi_t} \log  \left[  1 + \xi_t \left(\frac{Z_t-\mu_t}{\sigma_t} \right) \right].
            \end{equation}
Provided the chosen covariates appropriately explain the non-stationarity in the parameters of the GEV corresponding to $Z_t$, the standardized variable $\Tilde{Z}_t$ is expected to follow a standard Gumbel distribution. A conventional quantile plot for block maxima can be generated using the standardized variable $\Tilde{Z}_t$ compared to a standard Gumbel distribution. Anderson-Darling (A-D) and Mann-Kendall (M-K) tests are used to assess how well the data fit the distribution, and how well our covariates explain the non-stationarity within the data. Model diagnostics, including M-K and A-D test results can be found in Table \ref{tab:diag} and quantile plots to assess the fit against a stationary model can be found in Figure \ref{fig:qqplot}.


\subsection*{Assessing model complexities and improvements for multiday extremes.}
We test each non-stationary model against its stationary counterpart to determine whether increased model complexity statistically significantly improves model performance. We refer the reader to \cite{KIM} which provides a nice summary on appropriate model selection criteria under nonstationary assumptions. In line with the results of this investigation, we use the AIC to select the appropriate model out of our series of nested models; however, we remark that likelihood ratio and BIC were also in agreement with this choice. Nested nonstationary models were used to assess the statistical significance against model complexities. These models are described in the appendix \ref{append}. In agreement with previous literature on single-day precipitation extremes across eastern Australia \cite{saunders2017spatial}, we find a non-stationary model with latitude, longitude, distance from the coast, and SOI as the most representative of the non-stationarity observed in multiday precipitation extremes.
    
\subsection*{Estimating return-levels and return-times under non-stationary conditions.}\label{rl}
In applications, we are usually interested in how often the extreme quantiles occur with a specific return level rather than parameter estimates. Return level ($z_p$) refers to a value expected to be equaled or exceeded on average once every interval of time or return period $1/p$ with probability $p$, defined by
        \begin{equation}
            z_p =
            \begin{cases}
                \mu - \frac{\sigma}{\xi}(1-y_p^{-\xi}), & \text{for } \xi \neq 0,\\
                \mu - \sigma \log y_p, & \text{for } \xi = 0,\\
            \end{cases}
        \end{equation}
    where $y_p=-\log(1-p)$. 
 
The estimation of return levels becomes difficult when underlying distribution properties vary with covariates, violating stationarity assumptions. To cope with this issue, several measures have been proposed \cite{cheng2014non, cooley2012return,katz2002statistics}. To estimate the non-stationary $m$-year return level, we used the idea of the aggregated quantile of the order $p$, $\Bar{x}_p$, following \cite{markiewicz2020quantile} defined as a weighted average of the form:
        \begin{equation}
            \Bar{x}_p = \sum_{i=i}^m w_i (\hat{x}_p)_i,
        \end{equation}
    where $(\hat{x}_p)_i$ is the quantile that is calculated for the $i$-th year out of the $m$ years that are being considered, and $w_i$ is the weight assigned to each subsequent year equally, given by $w_i=1/m$. In this setting, $\Bar{x}_p$ is interpreted as the return level one expects to exceed in $m$ years with probability equal to 1.
    
\section{Results and Conclusions}

\subsection*{Multiday precipitation extreme model predictions for historical and future flooding events.}
We first examine the spatial and temporal accuracy of our model to represent multiday extreme rainfall for years of historical relevance. Then, we report on the return level expectations under varying ENSO conditions.

\subsubsection*{Benchmarked return levels to empirical rainfall data from historical flooding events.} During the summer of 2010-2011, extensive parts of Queensland, New South Wales and Victoria experienced one of the worst flooding events on record, according to economic and loss of life measures. Following in 2022, Queensland and New South Wales experienced 2-3 consecutive days of extreme rainfall, causing extensive property damage, delaying harvests and reducing grain quality in the growing regions of New South Wales and Victoria. A natural question is whether we can use our non-stationary models of consecutive extremes to inform predictions on expected returns of consecutive days of extreme precipitation. 

   \begin{figure}[h!]
	 \begin{center}
	 \hspace*{-0.8in}
            \includegraphics[width=1.2\textwidth]{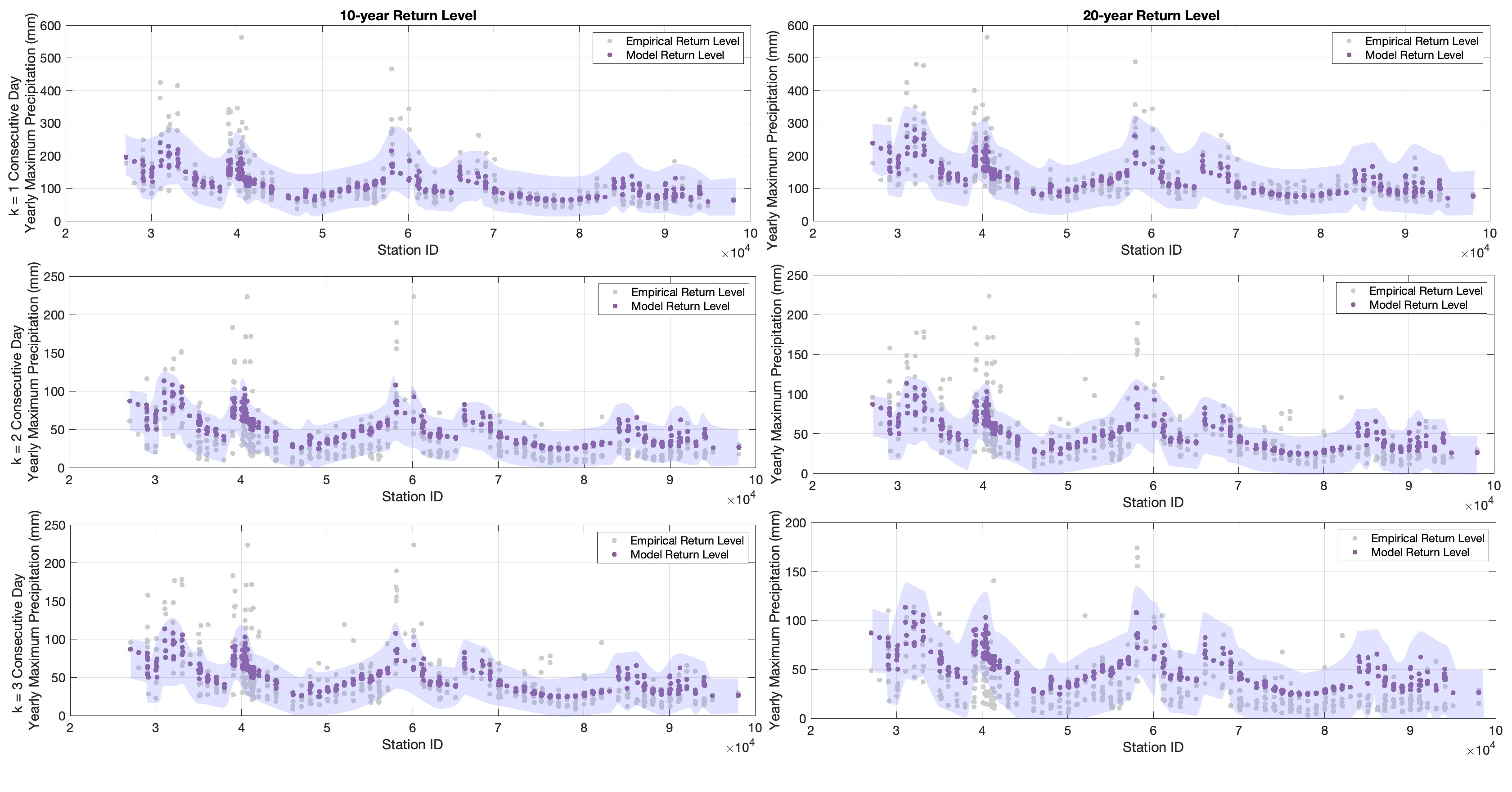}
            \end{center}
            \caption{2022 return levels, 10-years out from 2011 and 20-years out from 2001, for mutliday extreme rainfall across stations in eastern Australia for $k = 1$, $k = 2$, and $k = 3$ consecutive days (from top).}
            \label{fig:2022}
        \end{figure}

We adopt the approach in section \ref{rl} which allows us to answer the question: \textit{What level of extreme rainfall do we expect to observe for $k$ = 1, 2, or 3 consecutive days 10 years out (in 2022) from the last flooding event (in 2011)?} Return level estimates and 95\% confidence intervals are shown in Figure \ref{fig:2022} against empirical quantiles. In noting that these empirical quantile estimates are more accurate with increasing data, we also provide the 20 year return level for 2022 (e.g. predictions out from 2001). Our models agree well with the empirical quantiles of multiday extremes in 2022.

\subsubsection*{Model behaviour for historical El Ni\~{n}o and La Ni\~{n}a conditions.} 
We use these methods to illustrate differences in the probability of observing multiday extreme precipitation under historical El Ni\~{n}o and La Ni\~{n}a conditions. We note that our model results suggest the contribution of SOI on the variance for $k = 3$ consecutive days is double that of $k = 1$ or $2$ consecutive days of extreme precipitation, against an overall decreasing trend for variance against $k$, which agrees well with the conclusions of \cite{westra2011detection}, suggesting the SOI is most influential for longer duration storms. Our model also suggests that changes to SOI have greater effects on inland compared to coastal regions. Table \ref{tab:soieffect} shows the $k = 1$, $2$ and $3$ estimates of consecutive extreme rainfall using El Ni\~{n}o and La Ni\~{n}a conditions from 2015 and 2010, respectively. As an example, Karange Model Farm reports a 4.3\% increase in probability of observing 3 consecutive days of extreme rainfall above the regionally high threshold compared to approximately 1-1.5\% increase for coastal regions. The percent increase in observing multiday extreme rainfall from El Ni\~{n}o to La Ni\~{n}a years increases across duration (the number of consecutive days).

\begin{table}[h!]
\begin{center}
\begin{tabular}{ c | c c c c c c c c c c c}
 & Townsville & Brisbane & Sydney & Karange\\
\hline\hline
1-Day Consecutive Extreme Daily Rainfall (mm) & 211.09 & 178.70 & 155.15 & 66.04\\
La Ni\~{n}a Probability of Exceedance (2010) & 0.10 & 0.10 & 0.10 & 0.10\\
El Ni\~{n}o Probability of Exceedance (2015) & 0.09 & 0.09 & 0.09 & 0.08\\
Percent Increase in Probability & 0.74 & 0.93 & 1.06 & 2.5\\
\hline\hline
2-Day Consecutive Extreme Daily Rainfall (mm) & 99.60 & 88.44 & 76.52 & 26.3\\
La Ni\~{n}a Probability of Exceedance (2010) & 0.10 & 0.10 & 0.10 & 0.10\\
El Ni\~{n}o Probability of Exceedance (2015) & 0.09 & 0.09 & 0.09 & 0.07\\
Percent Increase in Probability & 0.85 & 1.01 & 1.17 & 3.23\\
\hline\hline
3-Day Consecutive Extreme Daily Rainfall (mm) & 54.60 & 50.71 & 44.84 & 12.74\\
La Ni\~{n}a Probability of Exceedance (2010) & 0.10 & 0.10 & 0.10 & 0.10\\
El Ni\~{n}o Probability of Exceedance (2015) & 0.09 & 0.09 & 0.09 & 0.06\\
Percent Increase in Probability & 1.08 & 1.20 & 1.37 & 4.30
\end{tabular}
\end{center}
\caption{Probabilities of multiday extreme rainfall exceedances for select stations and 2010, 2015 ENSO conditions.}\label{tab:soieffect} 
\end{table}

   \begin{wrapfigure}{l}{0.48\textwidth}
	 \begin{center}
            \includegraphics[width=0.48\textwidth]{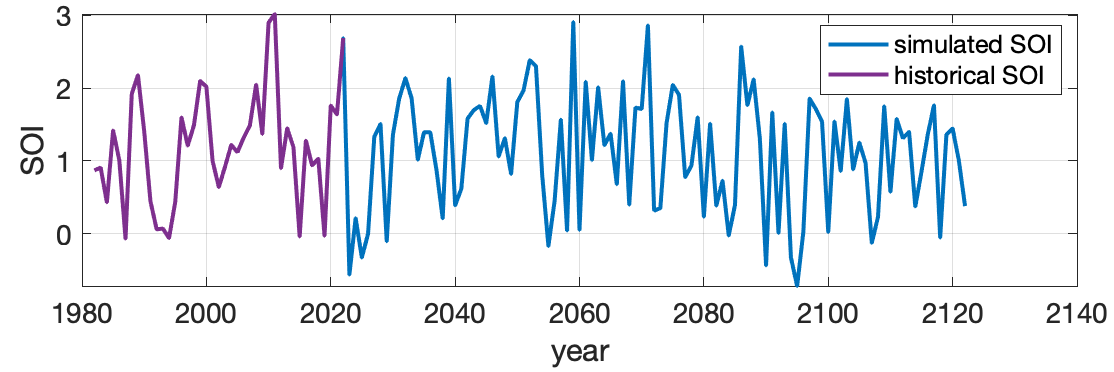}
             \vspace{-20pt}
            \end{center}
            \caption{Example simulated SOI from randomly sampling from the distribution of averaged yearly historical SOI.}
            \label{fig:soisim}
             \vspace{-10pt}
        \end{wrapfigure}

\subsubsection*{Predicted return levels under varying ENSO conditions.} 
\begin{figure}
	 \begin{center}
            \includegraphics[width=\textwidth]{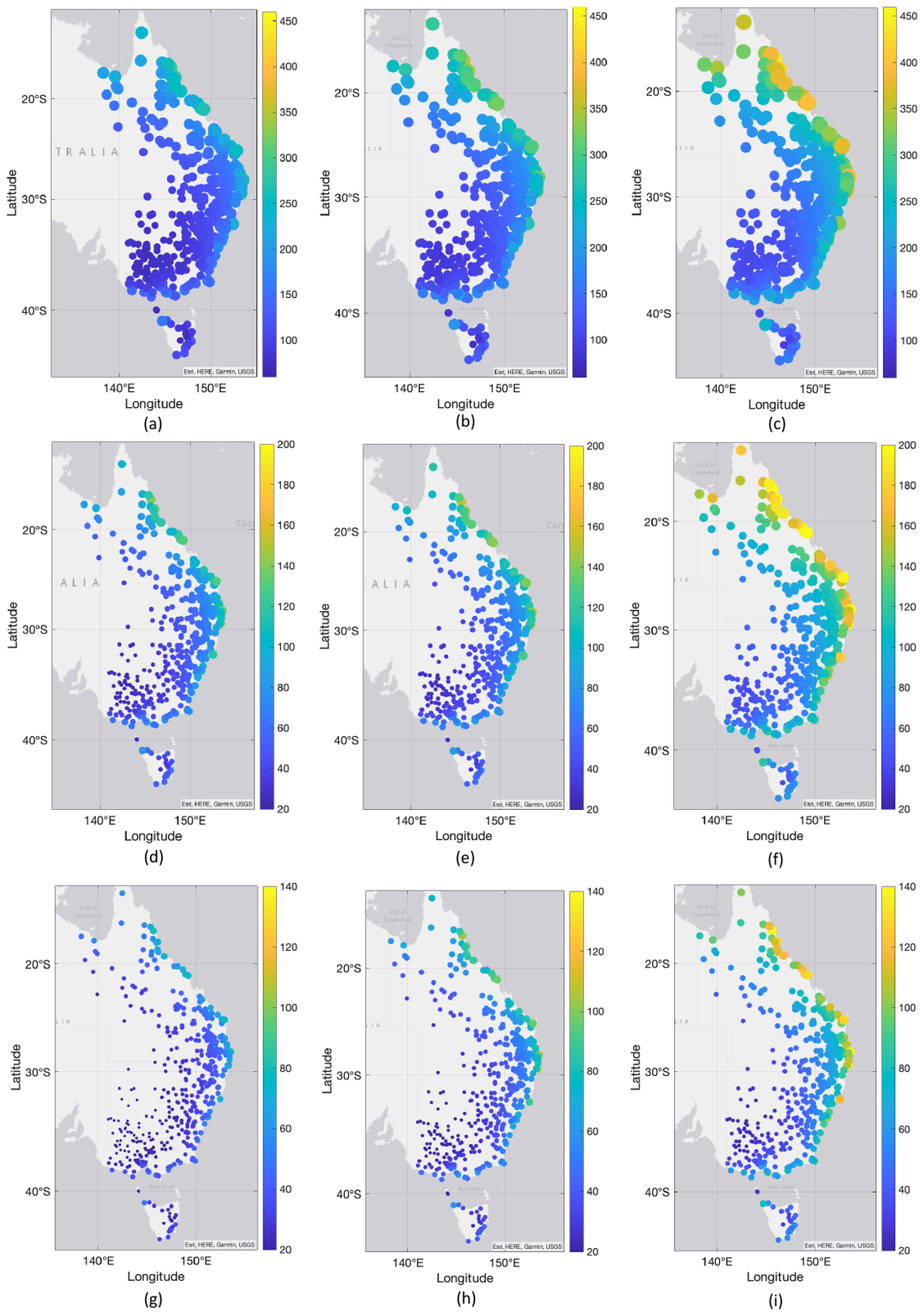}
            \end{center}
            \caption{From the left to right, 20-, 40-, and 100-year return levels for (top to bottom) $k=1$, 2, and 3 consecutive days of extreme rainfall.}
            \label{fig:final}
        \end{figure}

For the readers' interest, we simulate varying ENSO conditions for 20, 40, and 100 years by random sampling from the distribution of averaged yearly historical SOI taken from the NOAA (Figure \ref{fig:soisim}). We emphasize that these simulations are meant to provide ideas for tools and techniques for climate scientists and should be interpreted in that context. Return levels for multiday extreme rainfall from these simulations for all locations along the eastern coast of Australia can be found in Figure \ref{fig:final}.

\subsubsection*{A note on the prediction horizon.} The prediction accuracy of our non-stationary model for returns and magnitudes of multiday precipitation extremes across eastern Australia is dependent on the accuracy of predictions of the proceeding years' ENSO (SOI or similar) index.

\subsection*{Conclusions and further applications.}
We have summarized current, popular methods used to model extreme rainfall over varying durations and the relationship of these methods to the theoretical advancements made in the dynamical framework. We have shown where classical extreme value modeling fails to appropriately interpret returns of consecutive extreme events (e.g. through averaging over the extremal index). We have introduced the minimum functional, taken from dynamical systems, as a way to preserve the consecutive property of multiday extremes in rainfall. We illustrate how advancements in dynamics can be translated to improve hidden numerical error (e.g. not detectable by current statistical methods) on multiday extreme rainfall models.

 We have shown that in the case of extreme rainfall following a Fr\'{e}chet distribution $(\xi>0)$, the max-stability results from \cite[Lemma 4.1]{CHNT} can be used to improve model estimates for returns of multiday extreme rainfall by limiting investigations to window sizes which guarantee a stabilized shape parameter. This is the first time max-stability properties have been used to appropriately model \textit{consecutive} extreme rainfall in a non-stationary setting; however, we note that these results are complementary to those observed empirically in previous literature on extremes of rainfall taken over higher resolutions, but shorter durations where only small fluctuations in shape were observed \cite{westra2011detection}. We have shown that our models agree with the empirical historical estimates of rainfall and the conclusions of these models agree with current literature, in particular on the empirically observed effects of SOI on extreme rainfall duration. Although we have used eastern Australia as a case study, we note the general use of these methods for modeling multiday rainfall in many locations where the maxima of daily rainfall follow a Fr\'{e}chet distribution.

The importance of the results of this investigation is emphasized by inconsistent empirical results on shape parameter estimates for extreme rainfall over varying durations reported in recent literature. For example, the authors in \cite{westra2011detection} reported a positive shape parameter with some variation depending on storm duration (6 minutes to 72 hours) across Australia. The authors in \cite{HR} report that their empirically estimated shape parameters range from slightly negative to positive values for 1-3 days durations at sites in New South Wales. On the other hand, the authors in \cite{yilmaz2017investigation} obtained negative shape parameters for 1- and 2-day durations at a specific site in Victoria and \cite{JOHNSON201867} reports a negative shape parameter for 1-day durations is common across meteorological stations in Australia. It is not clear whether the negative shape parameters observed for these regions is due to a true upper bounds on the quantity of rainfall for the chosen durations or whether it is the result of lack of available data. We illustrate that the data for modelling extreme multiday rainfall follow the general assumptions of \cite[Lemma 4.1]{CHNT}, so that there is no theoretical expectation for a varying shape parameter over storm duration and that the observation of such a change may be the result of increased numerical error in limiting the available pool of data used to fit the GEV. 

The authors in \cite{CHNT} have also investigated max-stability for averaging functionals (e.g. extremes of the average over some window of length $k$). In future investigations, we would like to develop techniques for modeling heat-waves using such results. Due to the novelty of these results, convergence rates for appropriate block-size choice have not yet been investigated. We intend on performing numerical investigations to estimate these rates which can then be used to choose appropriate block-sizes against window size $k$ (duration of consecutive extreme events) and perhaps incorporated into more robust statistically fitting methods for consecutive extreme events. 

\clearpage
\appendix
\section{Supplementary figures and tables}\label{append}

        \begin{figure}[h]
            \centering
            \includegraphics[width=\textwidth]{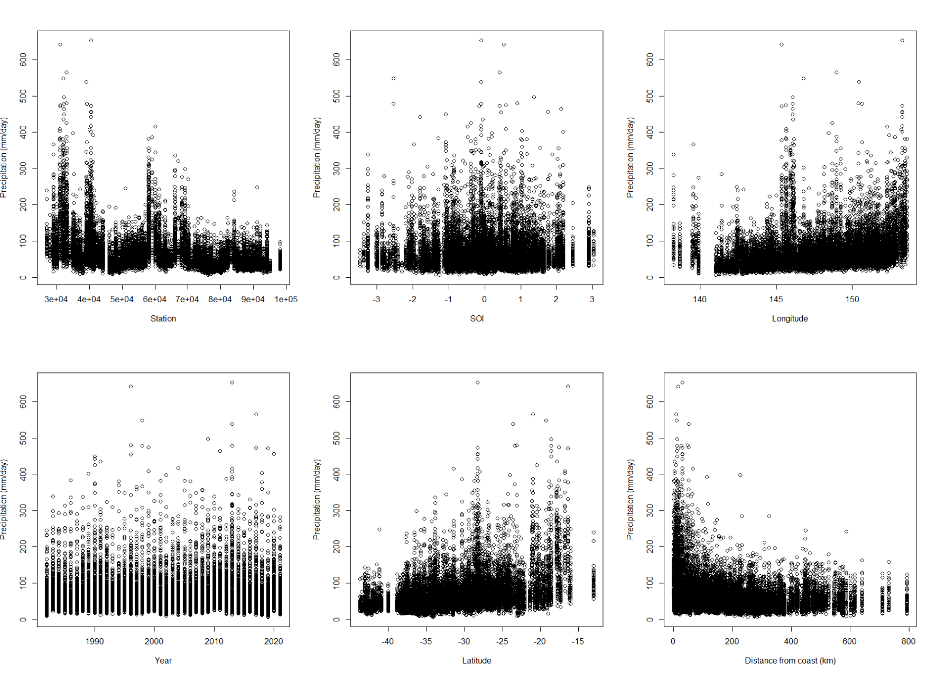}
            \caption{Scatter plots of the block maxima against covariates assess for the nonstationary models.}
            \label{fig:scatter}
        \end{figure}
        
         \begin{figure}[h]
            \centering
            \includegraphics[width=\textwidth]{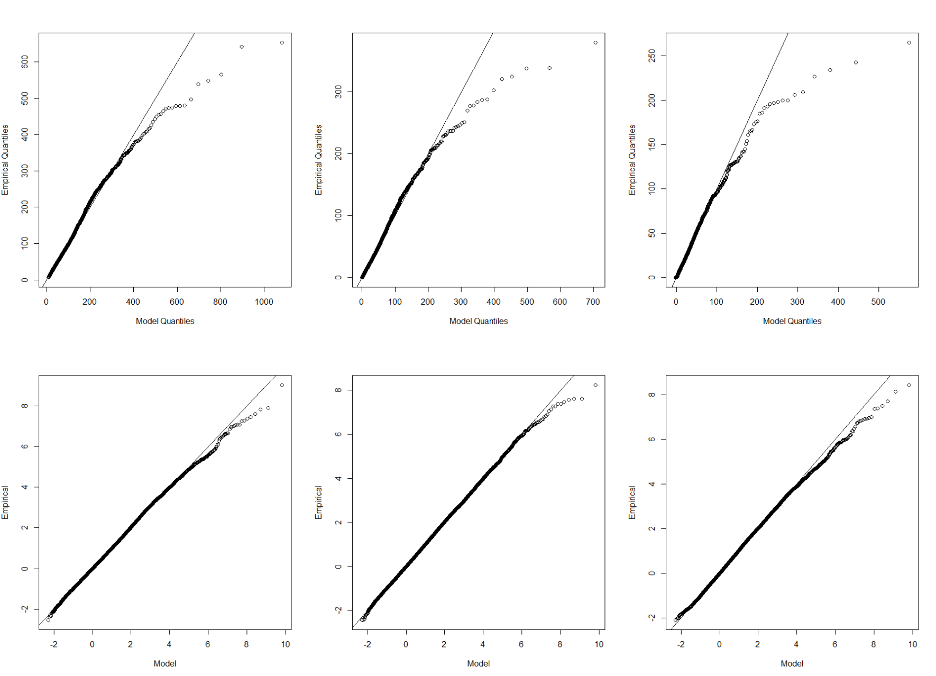}
            \caption{Quantile plots of the stationary (top) and chosen nonstationary (bottom) models.}
            \label{fig:qqplot}
        \end{figure}

\begin{table}[h!]
\begin{center}
\begin{tabular}{ c | c c c c c c c c c c c}
window size & $\mu_0$ & $\mu_1$& $\mu_2$ & $\mu_3$ & $\mu_4$ & $\sigma_0$ & $\sigma_1$ & $\sigma_2$ & $\sigma_3$ & $\sigma_4$ & $\xi$\\
\hline \hline
$k$ & int & SOI & $\log(\text{cdist})$ & lat & lon & int & SOI & $\log(\text{cdist})$ & lat & lon & int\\
\hline
1 &	8.40 &	0.12 &	0.16	 &	0.04	&	0.05	&	7.07 &	0.09	&	0.13	 &	0.03	&	0.05	&	0.01\\
2 &	3.91 &	0.05 &	0.08 &	0.02  &	0.02	&	3.36	&	0.04	&	0.07	 &	0.01	&	0.02	&	0.01\\
3 &	2.32 &	0.03 &	0.05 &	0.01	&	0.01	&	3.24	 &	0.02	&	0.04	 &	0.01	&	0.02	&	0.01\\
\end{tabular}
\end{center}
\caption{Chosen nonstationary model coefficient standard errors. \label{tab:se}}.
\end{table}

\begin{table}[h!]
\begin{center}
\begin{tabular}{ c | c c c c c c c c}
model type & window size & NLL & AIC & BIC & M-K & & A-D & \\
\hline \hline
&  & & & & statistic & p-value & statistic & p-value\\
\hline
0	&1	   &88,509.51 	  &177,025.00 	  &177,048.40 	&0.99971	&0.00000	&3.04977	&<0.001\\
1	&1	   &88,474.22 	  &176,958.40 	  &176,997.50 	&0.01088	&0.02660	&3.14282	&<0.001\\
2	&1	   &87,697.31 	  &175,404.60 	  &175,443.70 	&0.00397	&0.93663	&0.48007	&0.2398\\
3	&1	   &84,325.64 	  &168,669.30 	  &168,739.60 	&0.00851	&0.14232	&1.20606	&0.0038\\
4	&1	   &84,250.02 	  &168,522.00 	  &168,608.00 	&0.00684	&0.36143	&1.13899	&0.0057\\
\hline
0	&2	   &75,417.05 	  &150,840.10 	  &150,863.50 	&0.97821	&0.00000	&5.59706	&<0.001\\
1	&2	   &75,354.55 	  &150,719.10 	  &150,758.10 	&0.01302	&0.00411	&4.93463	&<0.001\\
2	&2	   &74,268.54 	  &148,547.10 	  &148,586.10 	&0.00474	&0.80760	&0.87351	&0.0275\\
3	&2	   &71,484.01 	  &142,986.00 	  &143,056.30 	&0.00521	&0.70548	&0.60843	&0.1202\\
4	&2	   &71,400.54 	  &142,823.10 	  &142,909.00 	&0.00519	&0.70923	&0.67106	&0.0854\\
\hline
0	&3	   &64,543.39 	  &129,092.80 	  &129,116.20 	&0.84080	&0.00000	&6.09299	&<0.001\\
1	&3	   &64,469.64 	  &128,949.30 	  &128,988.30 	&0.01065	&0.03198	&4.83873	&<0.001\\
2	&3	   &62,937.73 	  &125,885.50 	  &125,924.50 	&0.01579	&0.00022	&6.99238	&<0.001\\
3	&3	   &61,101.78 	  &122,221.60 	  &122,291.90 	&0.01007	&0.04949	&2.10993	&<0.001\\
4	&3	   &61,032.96 	  &122,087.90 	  &122,173.80 	&0.00955	&0.07157	&1.78479	&<0.001\\

\end{tabular}
\end{center}
\caption{Model diagnostic and model selection: log-likelihood from maximum likelihood estimation, AIC, BIC, Likelihood ratio test, and goodness of fit tests for nested nonstationary models. \label{tab:diag}}.
\end{table}
\clearpage
\section*{Declarations}
\subsection*{Ethical Approval}
The authors declare that no ethical approvals are required for this study.
\subsection*{Availability of supporting data}
All supporting data is freely available from the Australian Bureau of Meteorology website through the following link: \\http://www.bom.gov.au/climate/data/
\subsection*{Competing interests}
The authors have no competing interests to declare.
\subsection*{Funding}
The authors would like to acknowledge the University of Queensland for facilities and funding under the research start-up scheme.
\subsection*{Authors' contributions}
M.C. wrote the main manuscript text, formatted all tables and figures, prepared figures 3, 4, and 5, and completed the formal analysis. R.K. wrote the main manuscript text, prepared figures 1, 2, 6 and 7, prepared tables 1, 3, and 4, and completed the formal analysis. Both authors worked together on sections of the code (R and MATLAB) for the numerical results and reviewed the manuscript.
\subsection*{Acknowledgements}
M Carney would like to thank the Max Planck Institute for the Physics of Complex Systems, Dresden, Germany, where part of this work was completed.
\clearpage

\bibliographystyle{plain}
\bibliography{Reference}

\end{document}